\documentclass[11pt]{article}
\usepackage{amsmath}
\usepackage{amssymb}
\usepackage{amscd}
\usepackage{amsmath}
\usepackage[all]{xy}
\usepackage{setspace}
\usepackage{graphics,epsfig}
\usepackage{graphicx}
\usepackage{float}
\usepackage{stackrel}
\usepackage{tikz-cd}
\usepackage{tikz}
\usetikzlibrary{positioning}
\usepackage{amsthm}
\usepackage{stix}
\newcommand{\sixvdots}{%
  {\vbox{\baselineskip1ex\lineskiplimit0pt%
  \hbox{.}\hbox{.}\hbox{.}\hbox{.}\hbox{.}}}}

\newcommand{\tenvdots}{%
  \vbox{\baselineskip1ex\lineskiplimit0pt%
  \hbox{.}\hbox{.}\hbox{.}\hbox{.}\hbox{.}\hbox{.}\hbox{.}\hbox{.}\hbox{.}\hbox{.}\hbox{.}\hbox{.}\hbox{.}}}

\DeclareMathAlphabet{\eusm}{OT1}{eusm}{m}{n}
\newtheorem{thm}{Theorem}[section]

\newtheorem{exam}[thm]{Example}
\newtheorem{rem}[thm]{Remark}

\def\vsp{\vspace{1ex}}

\begin{document}
\begin{center}
{\rm {\LARGE Bijections between $\tau$-rigid modules}}
\end{center}
\
\

\begin{center}
{\rm {\Large Gabriella D$^{'}$Este$^a$ and H. Melis Tekin Akcin$^b$
}}\\
${}^a$ Department of Mathematics, University of Milano,\\
Milano, Italy\\
e-mail: gabriella.deste@unimi.it\\

${}^b$ Union Lane CB4 1QB, \\
Cambridge, United Kingdom\\
e-mail: meliss85@gmail.com\\
\end{center}

\vsp
\begin{abstract}

We describe a special bijection between the indecomposable summands of two basic $\tau$--tilting modules.

{\em Key Words:  Tilting modules, $\tau$-tilting modules, $\tau$-rigid modules, quivers.}

{\em Mathematics subject classification 2020: Primary: 16G20; Secondary: 16D10.}

\end{abstract}

\section{Introduction}

There are many types, more or less abstract, of tilting objects which generalize the tilting modules first studied by Brenner and Butler in \cite{BB}. In this note, we consider $\tau$-tilting modules in the sense of \cite{IR} and tilting modules in the sense of \cite{IR} and \cite[page 167]{R}. More precisely, let $K$ be a field and let $A$ be a finite dimensional $K$-algebra with $n$ simple modules. Following and \cite{AIR} and \cite{IR}, we say that a module $T$ is 
\begin{itemize}
\item[$\bullet$] \emph{$\tau$-rigid} if $Hom_{A}(T, \tau(T))=0$;

\item[$\bullet$] \emph{$\tau$-tilting} if $T$ is $\tau$-rigid and $T$ is the direct sum of $n$ pairwise non-isomorphic indecomposable modules.

Moreover, following \cite{IR} and \cite[page 167]{R}, we say that a module $T$ is 
\item[$\bullet$] \emph{partial tilting} if the projective dimension of $T$ is at most one and \\$Ext^{1}_{A}(T,T) = 0$. 

\item[$\bullet$] \emph{tilting} if $T$ is partial tilting and the number of the isomorphism classes of the indecomposable summands of $T$ is equal to $n$.
\end{itemize} 
We know from \cite[Proposition 10(a)]{IR} that a $\tau$-tilting module $T$ is faithful if and only if it is a \emph{basic} (or \emph{multiplicity free}) \emph{tilting} module, that is a tilting module in the above sense which is the direct sum of pairwise non-isomorphic indecomposable summands.
We know from \cite[Theorem 2.5]{GM} that, given two basic tilting modules $X= \oplus_{i=1}^{n} X_{i}$ and 
$Y=\oplus_{i=1}^{n} Y_{i}$ with $X_{i}$ and $Y_{i}$
indecomposable for any $i=1,\ldots, n$, there is a permutation $s$ of $\{1, \ldots, n\}$ with the following properties:
\begin{itemize}
\item[$\bullet$] If $X_{i}\simeq Y_{j}$, then we have $s(i)=j$.
\item[$\bullet$] If $X_{i}\not\simeq Y_{j}$ for any $j=1, \ldots, n$, then we have that \[ Ext^1_{A}(X_{i},Y_{s(i)})\oplus 
Ext^1_{A}(Y_{s(i)},X_{i})\neq 0.\]
\end{itemize}
We will show that we can not replace the hypothesis that both $X$ and $Y$ are basic tilting modules with the weaker hypothesis that they both are $\tau$-tilting modules (Theorem \ref{Theorem5}). We will also show that an algebra $A$ and its opposite may have the same number of $\tau$-tilting modules, but a different number of basic tilting modules (Examples \ref{Ex 2} and \ref{Ex 3}).

In the following, all modules are left modules defined over algebras given by quivers. Moreover, for any vertex $x$, we denote by $P_{x}$, $I_{x}$ and $S_{x}$ the indecomposable projective, injective and simple modules corresponding to $x$. 

\section{Examples and Results}
\begin{exam}\label{Ex 1}
There are a $K$-algebra $A$ and three non-faithful $\tau$-tilting $A$-modules $X$, $Y$ and 
$Z$ with the following properties:

\begin{itemize}
\item[(i)]If $U, V\in \{X,Y,Z\}$ and $U\not\simeq V$, then $U$ has
an indecomposable projective-injective summand which is not a summand of $V$. 
\item[(ii)] $X$, $Y$ and $Z$ are sincere and their annihilators are ideals of $A$ of dimension one over $K$.
\end{itemize}
\end{exam}
\noindent
{\bf Construction:} Let $A$ be the algebra given by the quiver

$$\begin{tikzcd}
\underset{1}{\bullet} \arrow[r, "a"] & \underset{2}{\bullet} \arrow[r,"b"] & \underset{3}{\bullet} \arrow[ll, "c", bend left=30]
\end{tikzcd}$$
with relations $ba=0$, $cb=0$ and $ac=0$. Then the Auslander-Reiten quiver of $A$ has the following form: 

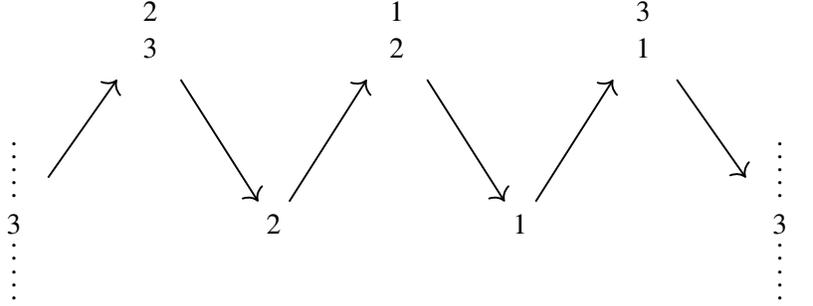
\begin{center}
$\begin{tikzcd}
   &\begin{array}{c} 2\\3 \end{array} \arrow[rd] &  &\begin{array}{c} 1\\2 \end{array} \arrow[rd] & & \begin{array}{c} 3\\1 \end{array} \arrow[rd] &  \\
\begin{array}{c}\sixvdots \\ 3 \\ \sixvdots \end{array} \arrow[ru] &  & 2 \arrow[ru] & & 1 \arrow[ru]  & & \begin{array}{c}\sixvdots \\ 3 \\ \sixvdots \end{array} 
\end{tikzcd}.$
\end{center}
It follows that any $\tau$-tilting module has at most one simple summand. This observation implies that \\
%$$Ext^1_{A}(1,1)= 0, Ext^2_{A}(1,1) \simeq Ext^1_{A}(2,1)=0,$$ 
%$$Ext^3_{A}(1,1)\simeq Ext^2_{A}(2,1)\simeq Ext^1_{A}(3,1)\neq 0.$$
$(1)$ The non-faithful $\tau$-tilting modules are $X= \begin{matrix} 1\\2\end{matrix}\oplus \begin{matrix}
2\\3 \end{matrix}\oplus 2$, $Y = \begin{matrix} 2\\3 \end{matrix}\oplus \begin{matrix}
3\\1 \end{matrix}\oplus 3$, $Z= \begin{matrix} 3\\1 \end{matrix}\oplus \begin{matrix} 1\\2 \end{matrix}\oplus 1.$ \\  
Hence $(i)$ holds. Moreover, the modules $X, Y, Z$ are sincere and their annihilators are the one dimensional vector 
spaces generated by $c, a$ and $b$, respectively.
Hence, also $(ii)$ holds.\\

\begin{exam}\label{Ex 2} 
There exist an algebra $A$ and $A$-modules $X, Y$ and $Z$ with the following properties: 
\begin{itemize}
\item[(i)] $X$ is projective, while $Y$ and $Z$ have infinite projective dimension.
\item[(ii)] $X, Y$ and $Z$ are the $\tau$-tilting $A$-modules. Moreover $X$ is faithful, while $Y$ and $Z$ are not faithful.
\item[(iii)] There exists (resp. does not exist) a bijection $s$ between the indecomposable summands of $X$ and the indecomposable summands of $Y$(resp. $Z$), such that for any indecomposable summand $X_{i}$ of $X$, either $s(X_{i})$ is isomorphic to $X_{i}$ or $Ext^1_{A}(X_{i},s(X_{i}))\oplus Ext^1_{A}(s(X_{i}),X_{i})\neq 0$.
\item[(iv)] There is a permutation $t$ of the indecomposable $\tau$-rigid $A$-modules such that if $U$ and $W$ are $\tau$-tilting modules, then we obtain $W$ from $U$ by replacing any indecomposable summand $U_{i}$ of $U$ which is not a summand of $W$ by $t(U_{i})$.
\end{itemize}
\end{exam}
\noindent
{\bf Construction:} Let $A$ be the algebra given by the quiver 
$$\begin{tikzcd}
\underset{1}{\bullet} \arrow[r, "a"] & \underset{2}{\bullet} \arrow[r,"b"] & \underset{3}{\bullet} \arrow[out=120, in=50, loop right, "c", distance=1cm]
\end{tikzcd}$$
with relations $xy=0$ for any $x, y \in \{a, b, c\}$. Then the Auslander-Reiten quiver of $A$ is of the following form:

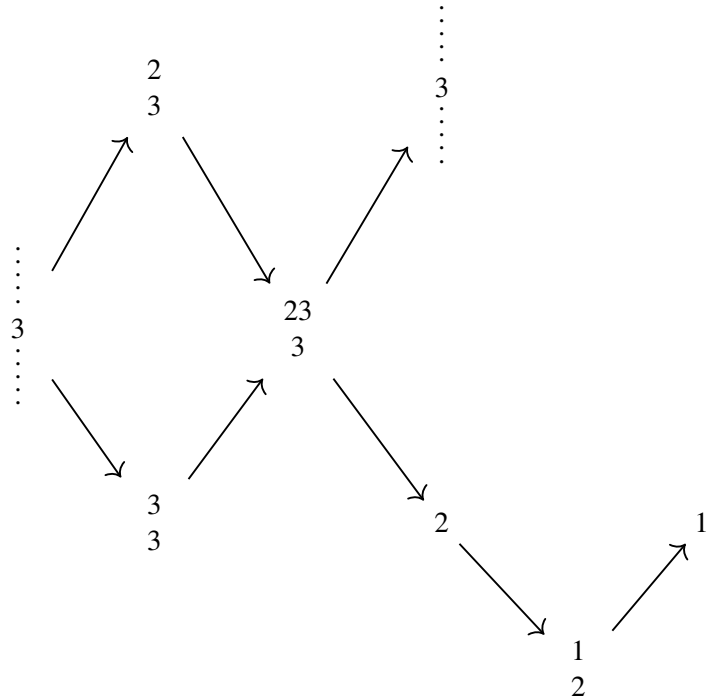
\begin{center}
$\begin{tikzcd}
   &\begin{array}{c} 2\\3 \end{array} \arrow[rd] &  & \begin{array}{c}\sixvdots \\ 3 \\ \sixvdots \end{array} & &  \\
\begin{array}{c}\sixvdots \\ 3 \\ \sixvdots \end{array} \arrow[ru] \arrow[rd] &  & \begin{array}{c} 23\\3 \end{array} \arrow[rd] \arrow[ru] & & & & \\
 & \begin{array}{c} 3\\3 \arrow[ru] \end{array} & & 2 \arrow[rd] & & 1 \\
& & & & \begin{array}{c} 1\\2 \arrow[ru] \end{array} & 
\end{tikzcd}.$
\end{center} 
The following facts hold: \\
$(1)$ $\begin{array}{c} 1\\2 \end{array}, \begin{array}{c} 2\\3 \end{array}, \begin{array}{c} 3\\3 \end{array}, 1$ and $2$
are the indecomposable $\tau$-rigid modules.\\
$(2)$ $1\oplus 2$, $1\oplus \begin{array}{c} 2\\3 \end{array}$, $2\oplus \begin{array}{c} 3\\3 \end{array}$ are not 
$\tau$-rigid.\\
$(3)$ The following modules are $\tau$-tilting:\\
$X=\begin{array}{c} 1\\2 \end{array} \oplus \begin{array}{c} 2\\3 \end{array} \oplus \begin{array}{c} 3\\3 \end{array}$,
$Y= \begin{array}{c} 1\\2 \end{array} \oplus \begin{array}{c} 2\\3 \end{array}\oplus 2$, $Z= \begin{array}{c} 1\\2 \end{array} \oplus \begin{array}{c} 3\\3 \end{array} \oplus 1$.\\ Then $X$ is projective, while the non-projective $\tau$-rigid modules $1$ and $2$ have infinite projective dimension. Hence $(i)$ holds. Let $T$ be a non-projective $\tau$-tilting module. Then we deduce from $(1)$ that one of the simple modules $1$ and $2$ is a summand of $T$. Assume first $T= 1\oplus \ldots.$ Then we deduce from $(1)$ and $(2)$ that $T= 1\oplus \begin{array}{c} 1\\2 \end{array}\oplus \begin{array}{c} 3\\3 \end{array}=Z$.
Suppose now $T= 2\oplus \ldots.$ Then we deduce from $(1)$ and $(2)$ that $T=2\oplus \begin{array}{c} 1\\2 \end{array} \oplus \begin{array}{c} 2\\3 \end{array}=Y.$ Hence $(ii)$ holds.   
Since $Ext_{A}^1\left(2, \begin{array}{c} 3\\3 \end{array}\right)\neq 0$ and $\begin{array}{c} 1\\2 \end{array} \oplus \begin{array}{c}2\\3 \end{array}$ is a common summand of $X$ and $Y$, it follows that \\
\\
$(4)$ The map $\begin{array}{c} 1\\2 \end{array} \mapsto \begin{array}{c} 1\\2 \end{array}$, $\begin{array}{c} 2\\3 \end{array} \mapsto \begin{array}{c} 2\\3 \end{array}$, $\begin{array}{c} 3\\3 \end{array} \mapsto 2$ is a bijection $s$ satisfying $(iii)$ between the indecomposable summands of $X$ and the indecomposable summands of $Y$.\\
\\
On the other hand, $\begin{array}{c} 1\\2 \end{array} \oplus \begin{array}{c} 3\\3 \end{array}$ is a common summand of $X$ and $Z$. We also note that any non-zero morphism $\begin{array}{c} 2\\3 \end{array}\rightarrow 2$ has a factorization of the form 
$
\begin{array}{c} 2\\3 \end{array} \hookrightarrow \begin{array}{c} 23\\3 \end{array}\twoheadrightarrow 2
$ 
with $\begin{array}{c} 23\\3 \end{array}$ injective. Since $2=\tau(1)$, we deduce from \cite[(5) page 75]{R} that 
$Ext_{A}^1\left(1, \begin{array}{c} 2\\3 \end{array}\right)=0.$ It follows that $Ext_{A}^1\left(\begin{array}{c} 2\\3 \end{array}, 1\right)\oplus 
Ext_{A}^1\left(1, \begin{array}{c} 2\\3 \end{array}\right)=0$. This implies that \\
\\
$(5)$ There is no bijection $s$ satisfying $(iii)$ between the indecomposable summands of $X$ and the indecomposable summands of $Z$.\\  
\\
Putting $(4)$ and $(5)$ together, we conclude that $(iii)$ holds. Finally, let $t$ denotes the following permutation of the indecomposable $\tau$-rigid $A$-modules:
$$ \begin{array}{c} 1\\2 \end{array} \longmapsto \begin{array}{c} 1\\2 \end{array}, \begin{array}{c} 2\\3 \end{array} \longmapsto \begin{array}{c} 1 \end{array}, \begin{array}{c} 3\\3 \end{array} \longmapsto \begin{array}{c} 2 \end{array}, \begin{array}{c} 1 \end{array} \longmapsto \begin{array}{c} 2\\3 \end{array}, \begin{array}{c} 2 \end{array} \longmapsto \begin{array}{c} 3\\3 \end{array}. $$ 
Then $t$ satisfies $(iv)$. \\

\begin{exam}\label{Ex 3} 
There are an algebra $B$ and $B$-modules $X,Y,Z$ with the following properties:\\
\begin{itemize}
\item[(i)] The projective dimension of $X, Y$ and $Z$ is $0,1$ and $2$, respectively.
\item[(ii)] $X, Y$ and $Z$ are the $\tau$-tilting $B$-modules. 
\item[(iii)] $X$ and $Y$ have an indecomposable projective-injective summand $P$ which is not a summand of $Z$. 
\item[(iv)] If $T\in \{X, Y\}$, then there is no bijection $s$ between the indecomposable summands of $T$ and 
the indecomposable summands of $Z$ such that, for any indecomposable summand $T_{i}$ of $T$, either $s(T_{i})$ is isomorphic 
to $T_{i}$ or $Ext_{B}^1(T_{i}, s(T_{i}))\oplus Ext_{B}^1(s(T_{i}), T_{i}))\neq 0$. 
\item[(v)] There is a permutation $t$ of the indecomposable $\tau$-rigid $B$-modules such that if $U$ and $W$ are $\tau$-tilting modules, then we obtain $W$ from $U$ by replacing any indecomposable summand $U_{i}$ of $U$ which is not a summand of $W$ 
by $t(U_{i})$.
\end{itemize}
\end{exam}
\noindent
{\bf Construction:} Let $B$ be the algebra given by the quiver 
$$\begin{tikzcd}
\underset{1}{\bullet} \arrow[r, "b"] \arrow[out=120, in=50, loop left, "a", distance=1cm] & \underset{2}{\bullet} \arrow[r,"c"] & \underset{3}{\bullet} 
\end{tikzcd}$$
with relations $xy=0$ for any $x, y \in \{a, b, c\}$. Then the Auslander-Reiten quiver of $B$ is of the following form:

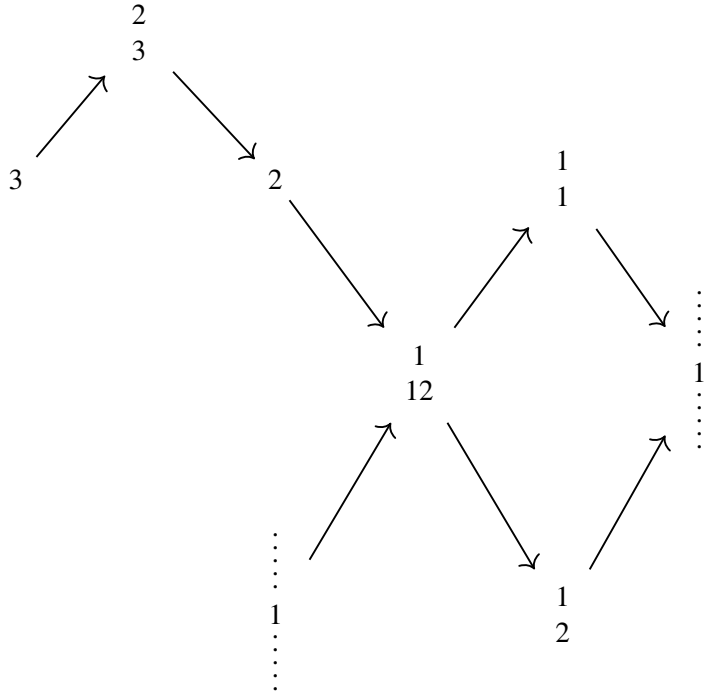
\begin{center}
$\begin{tikzcd}
  &\begin{array}{c} 2\\3 \end{array} \arrow[rd] &  &  & & & \\
3 \arrow[ru] &  & 2 \arrow[rd] & & \begin{array}{c} 1\\1 \end{array} \arrow[rd]&  \\
& & & \begin{array}{c} 1\\12  \end{array} \arrow[ru] \arrow[rd] &  & \begin{array}{c}\sixvdots \\ 1 \\ \sixvdots \end{array} \\
& & \begin{array}{c}\sixvdots \\ 1 \\ \sixvdots \end{array} \arrow[ru]& & \begin{array}{c} 1\\2 \end{array} \arrow[ru] & 
\end{tikzcd}.$
\end{center} 
The following facts hold:\\
$(1)$ $\begin{array}{c} 1\\12  \end{array}, \begin{array}{c} 2\\3  \end{array}, 3, 2, \begin{array}{c} 1\\1  \end{array}$
are the indecomposable $\tau-$rigid modules. \\
$(2)$ $\begin{array}{c} 1\\1  \end{array}\oplus 2$, $\begin{array}{c} 1\\1  \end{array}\oplus \begin{array}{c} 2\\3  \end{array}$, $2\oplus 3$ are not $\tau$-rigid.\\
$(3)$ The following modules are $\tau$-tilting:
$$
X=\begin{array}{c} 1\\12  \end{array}\oplus \begin{array}{c} 2\\3  \end{array}\oplus 3, \ Y= \begin{array}{c} 1\\12  \end{array}\oplus \begin{array}{c} 2\\3  \end{array}\oplus 2, \  Z=\begin{array}{c} 1\\12  \end{array}\oplus 3 \oplus \begin{array}{c} 1\\1  \end{array}. 
$$ 
It immediately follows that $(i)$ holds. Let $T$ be a non-projective $\tau$-tilting module. Then we deduce from $(1)$ and $(2)$
that $T=\begin{array}{c} 1\\12  \end{array}\oplus \ldots.$ Assume first \\$T= \begin{array}{c} 1\\12  \end{array}\oplus 2\oplus \ldots .$ Then we deduce from $(2)$ that $T= \begin{array}{c} 1\\12  \end{array} \oplus 2 \oplus \begin{array}{c} 2\\3  \end{array}=Y.$ Suppose now $T= \begin{array}{c} 1\\12  \end{array} \oplus \begin{array}{c} 1\\1  \end{array}\oplus \ldots.$
In this case, we deduce from $(2)$ that $T= \begin{array}{c} 1\\12  \end{array} \oplus \begin{array}{c} 1\\1  \end{array} \oplus 3=Z.$ Since $T$ is not projective, we know from $(1)$ that either  $2$ or $\begin{array}{c} 1\\1  \end{array}$ is a summand of $T$. Consequently, $(ii)$ holds. Moreover, let $P$ denote the module $\begin{array}{c} 2\\3  \end{array}$. Then we deduce from $(3)$ and $(ii)$ that $(iii)$ holds. Since $P$ is an indecomposable projective-injective module, we have $Ext_{B}^1(P, M)\oplus Ext_{B}^1(M, P)= 0$ for any module $M$. On the other hand, $P$ is direct summand of the tilting modules $X$ and $Y$, but not of $Z$. Hence, also $(iv)$ holds. Finally, let $t$ be the following permutation of the indecomposable $\tau$-rigid $B$-modules:
$$ \begin{array}{c} 1\\12  \end{array} \longmapsto \begin{array}{c} 1\\12  \end{array}, \begin{array}{c} 2\\3  \end{array}\longmapsto \begin{array}{c} 1\\1  \end{array}, \begin{array}{c} 3 \end{array} \longmapsto \begin{array}{c} 2  \end{array}, \begin{array}{c} 1\\1  \end{array} \longmapsto \begin{array}{c} 2\\3  \end{array}, \begin{array}{c} 2  \end{array} \longmapsto \begin{array}{c} 3  \end{array}.$$ 
Then $t$ satisfies $(v).$

\begin{rem} \label{Remark 4}
\textnormal{Note that the five indecomposable $\tau$-rigid modules in Examples \ref{Ex 2} and \ref{Ex 3} are
$P_1, P_2, P_3, I_1, S_2.$ We also note that, in both cases, we have $X=P_1\oplus P_2\oplus P_3$, $Y=P_1\oplus P_2 \oplus S_2$,
$Z= P_1 \oplus P_3 \oplus I_1$ and $I_1$ is isomorphic to $P_1/S_2$.
Moreover, in both cases, the bijection $t$ (satisfying condition $(iv)$ of Example \ref{Ex 2} and condition $(v)$ of Example \ref{Ex 3}) acts as follows:
$$ P_1 \longmapsto P_1, P_2\longmapsto I_1, P_3\longmapsto S_2, I_1\longmapsto P_2, S_2 \longmapsto P_3.$$
} 
\end{rem}

\begin{thm}\label{Theorem5}
Let $X= \oplus_{i=1}^{n} X_{i}$, $Y= \oplus_{i=1}^{n} Y_{i}$ be $\tau$-tilting modules over an algebra $A$
and let $(*)$ denote the following property: \\

$(*)$ There is a permutation $s$ of $\{1,\ldots, n\}$ such that, for any $i$, we have either $X_{i}\simeq Y_{s(i)}$
or $Ext^1_{A}(X_{i}, Y_{s(i)}) \oplus Ext^1_{A}(Y_{s(i)}, X_{i}) \neq 0$.\\ 
\newline
The following facts hold:
\begin{itemize}
\item[(i)] If $X$ and $Y$ are faithful modules, that is tilting modules, then $(*)$ holds.
\item[(ii)] If $(*)$ holds and $X$ is a faithful module, then $Y$ is not necessarily faithful.
%\item[(c)] If $X$ and $Y$ are two non-isomorphic $\tau$-tilting $A$-modules, then $(*)$ does not hold.
\item[(iii)] There is an algebra $A$ such that if $(*)$ holds, then $X$ and $Y$ are both faithful.
\end{itemize} 
\end{thm}
\noindent

\begin{proof}
\begin{itemize}
\item[(i)] This follows from \cite[Theorem 2.5]{GM} and from the fact \cite[Proposition 10(a)]{IR} that 
faithful $\tau$-tilting modules are exactly basic tilting modules.
\item[(ii)] This follows from Example \ref{Ex 2}(iii), where $X$ denotes the module \\$\begin{array}{c} 1\\2 \end{array} \oplus \begin{array}{c} 2\\3 \end{array} \oplus \begin{array}{c} 3\\3 \end{array}$ and $Y$ denotes the module $\begin{array}{c} 1\\2 \end{array} \oplus \begin{array}{c} 2\\3 \end{array} \oplus 2$.
%\item[(c)] We first note that the algebra $A$ in Example \ref{Ex 1} admits exactly one faithful $\tau-$tilting module, namely the module
%$\begin{array}{c} 1\\2 \end{array} \oplus \begin{array}{c} 2\\3 \end{array} \oplus \begin{array}{c} 3\\1 \end{array}$.
%On the other hand, the non-faithful $\tau$-tilting modules are 
%$\begin{array}{c} 1\\2 \end{array} \oplus \begin{array}{c} 2\\3 \end{array} \oplus 2$, $\begin{array}{c} 2\\3 \end{array} \oplus \begin{array}{c} 3\\1 \end{array} \oplus 3$, $\begin{array}{c} 3\\1 \end{array} \oplus \begin{array}{c} 1\\2 \end{array} \oplus 1$. This implies that, if $X$
%and $Y$ are non-isomorphic $\tau$-tilting modules, then $X$ has an indecomposable projective-injective summand $P$ which is not a summand of $Y$. Hence $X$ and $Y$ do not satisfy $(*)$.
\item[(iii)] This follows from Example \ref{Ex 3}(iv), where the unique non-faithful $\tau$-tilting module is the module $\begin{array}{c} 1\\12 \end{array} \oplus  3 \oplus \begin{array}{c} 1\\1 \end{array}$, while the other $\tau$-tilting modules, that is $\begin{array}{c} 1\\12 \end{array} \oplus \begin{array}{c} 2\\3 \end{array} \oplus 3 $ and
$\begin{array}{c} 1\\12 \end{array} \oplus \begin{array}{c} 2\\3 \end{array} \oplus 2$, admit the projective-injective module $\begin{array}{c} 2\\3 \end{array}$ as a summand. 
\end{itemize}
\end{proof} 
The algebras considered up to now admit at least two non projective indecomposable $\tau$-rigid modules and at least three
$\tau$-tilting modules.
In the next example, we construct a family of algebras with only one non projective indecomposable $\tau$-rigid
module and only one non projective $\tau$-tilting module.

\begin{exam}\label{Ex 6} 
For any $n\geq 2$ there exist two algebras $A$ and $B\simeq A^{op}$ with the following properties:\\
\begin{itemize}
\item[(i)] $A$ (resp. $B$) admit $2n+1$ indecomposable modules, one indecomposable non projective 
$\tau$-rigid injective module $I$ and one non projective $\tau$-tilting module $T$ which is not (resp. is) a tilting module. 
\item[(ii)] The Auslander-Reiten sequence of $A$-modules (resp., $B$-modules) ending in $I$ is of the form 
$$
0 \rightarrow \tau(I) \rightarrow X \rightarrow I\rightarrow 0
$$
with $\tau(I)$ projective and $X$ injective (resp., projective).  
\end{itemize}
\end{exam}
\noindent
{\bf Construction:} Let $A$ be the algebra given by the quiver
$$\begin{tikzcd}
\underset{1}{\bullet} \arrow[r, "a"] & \underset{2}{\bullet} \arrow[out=120, in=50, loop right, "b", distance=1cm]
\end{tikzcd}$$
with relations $ba=0$ and $b^{n}=0$. Next, let $B$ be the algebra given by the quiver 
$$\begin{tikzcd}
\underset{1}{\bullet} \arrow[r, "a"] \arrow[out=120, in=50, loop left, "b", distance=1cm] & \underset{2}{\bullet} 
\end{tikzcd}$$
with relations $b^{n}=0$ and $ab=0$.
Then there are $2n+1$ indecomposable $A$-modules(resp. $B$-modules) of the following 
form
$$ 
1, \begin{array}{c} 1\\2  \end{array}, \begin{array}{c} 12\\2  \end{array},
\begin{matrix}
 & & \begin{array}{c} 2  \end{array} \\ 
 & \reflectbox{$\ddots$} & \\
 \begin{array}{c} 12\\2  \end{array} & &
\end{matrix}, 2, \ldots, \begin{array}{c} 2\\ \vdots \\2  \end{array} 
$$
(resp. $1, \begin{array}{c} 1\\1  \end{array},
\ldots, \begin{array}{c} 1\\ \vdots \\1  \end{array}, 2, \begin{array}{c} 1\\2 \end{array},
\begin{array}{c} 1\\12  \end{array}, \ldots,
\begin{matrix}
 & & \begin{array}{c} 1\\12  \end{array} \\ 
 & \reflectbox{$\ddots$} & \\
 1 & &
\end{matrix}$).

To find the indecomposable $\tau$-rigid modules, it suffices to draw the Auslander-Reiten quivers 
of $A$ and $B$. As an example, assume $n=3$.
The Auslander-Reiten quiver of $A$ has the following form: \\

%\begin{center}
%$\begin{tikzcd}
%& & \begin{array}{c} 2\\2\\2\\2 \end{array} \arrow[rd] & & 1 &\\
%\ldots  & \begin{array}{c} \vdots \\  2\\2\\2 \\ \tenvdots \end{array} \arrow[rd] \arrow[ru] &  & \begin{array}{c} ~~~~2\\ 1~2~~\\~2~ \end{array} \arrow[rd]\arrow[ru] &  &\\
%\begin{array}{c} 2\\2 \end{array} \arrow[rd]\arrow[ru] & \tenvdots & \begin{array}{c} ~~~~2\\1~~2\\~2~ \end{array} \arrow[rd] \arrow[ru] &  & \begin{array}{c} \vdots\\ 2\\2\\2 \\ \tenvdots \end{array} & \cdots \\
%\ldots & \begin{array}{c} \tenvdots\\ 12\\2 \\ \vdots \end{array} \arrow[ru] \arrow[rd] & & \begin{array}{c} 2\\2 \end{array} \arrow[ru] \arrow[rd] & \tenvdots & \cdots \\
%\begin{array}{c} 1\\2 \end{array} \arrow[ru] & & 2 \arrow[ru] \arrow[rd] & & \begin{array}{c} \sixvdots \\ 12\\2 \\ \sixvdots\\ \end{array} & \cdots \\
%& && \begin{array}{c} \sixvdots \\ 1\\2  \end{array} \arrow[ru] & &
%\end{tikzcd}$.
%\end{center} 

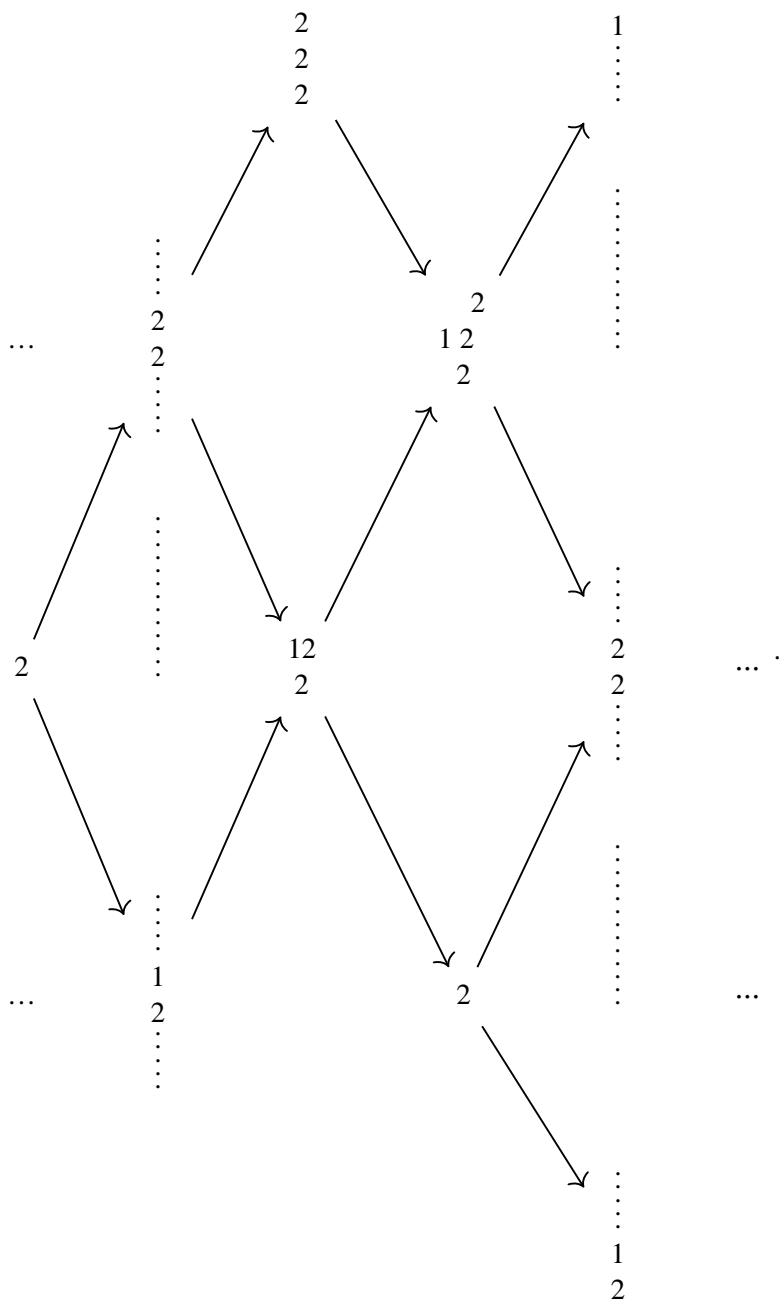
\begin{center}
$\begin{tikzcd}
& & \begin{array}{c} 2\\2\\2 \end{array} \arrow[rd] & & \begin{array}{c} 1 \\ \sixvdots \end{array}   &\\
\ldots  & \begin{array}{c} \sixvdots \\  2\\2 \\ \sixvdots \end{array} \arrow[rd] \arrow[ru] &  & \begin{array}{c}~~~~2\\ 1~2~~\\~2~ \end{array} \arrow[rd]\arrow[ru] & \tenvdots &\\
\begin{array}{c} 2 \end{array} \arrow[rd]\arrow[ru] & \tenvdots & \begin{array}{c} 12\\2 \end{array} \arrow[rd] \arrow[ru] &  & \begin{array}{c} \sixvdots\\ 2\\2 \\ \sixvdots \end{array} & \cdots \\
\ldots & \begin{array}{c} \sixvdots\\ 1\\2 \\ \sixvdots \end{array} \arrow[ru] & & \begin{array}{c} 2 \end{array} \arrow[ru] \arrow[rd]& \tenvdots & \cdots\\
& &  & & \begin{array}{c} \sixvdots \\ 1\\2 \end{array} &
\end{tikzcd}$.
\end{center} 

Consequently, the following facts hold: \\
$(1)$ The indecomposable $\tau$-rigid $A$-modules are the projective 
modules $\begin{array}{c} 1\\2 \end{array}$ and $\begin{array}{c} 2\\ \vdots \\2 \end{array}$
and the injective module $1$ of dimension $2$, $n$ and $1$ respectively. \\
$(2)$ The unique non projective $\tau$-tilting $A$-module is the non faithful module \\$T = \begin{array}{c} 1\\2 \end{array} \oplus 1$. \\
On the other hand, if $n=3$, then the Auslander-Reiten quiver of $B$ has the following form. 

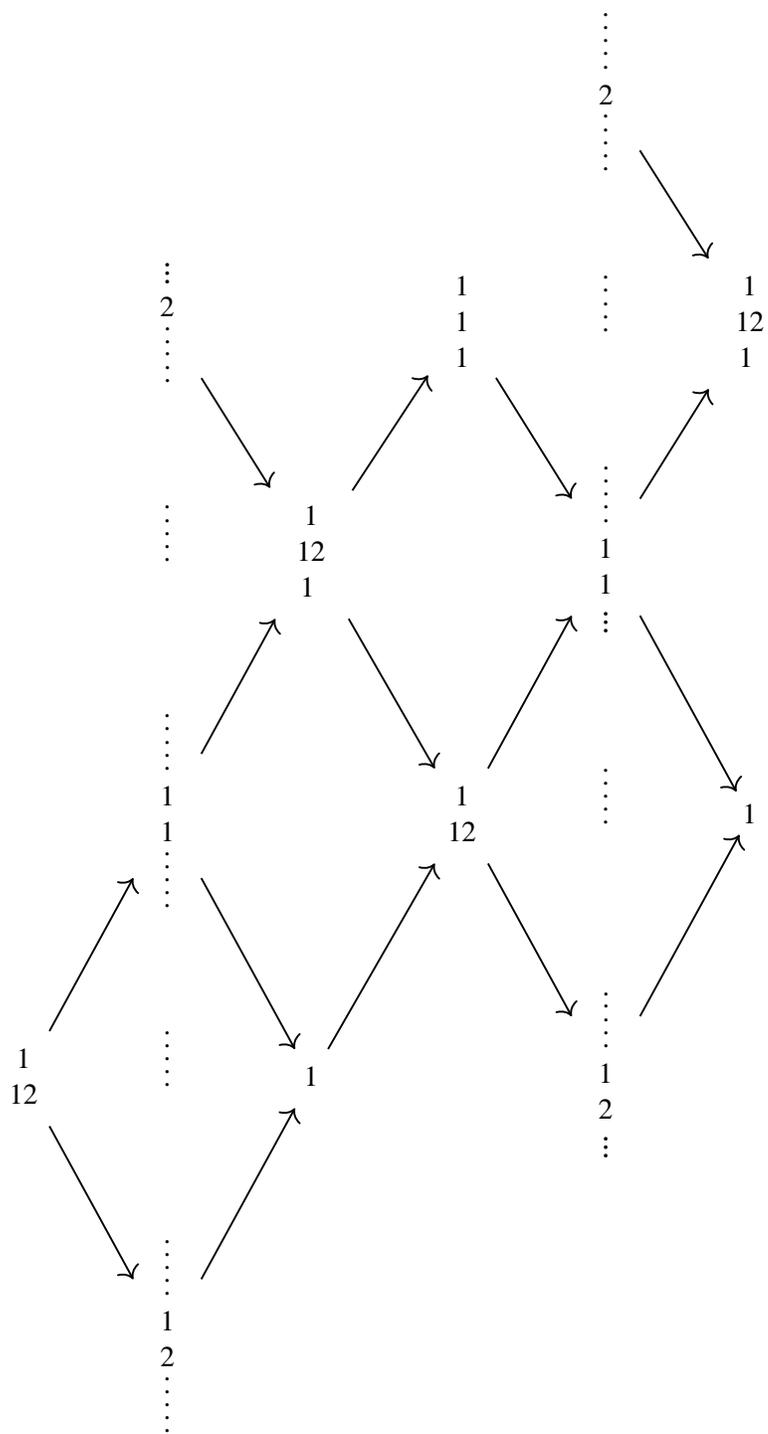
\begin{center}
$\begin{tikzcd}
& & &  & & & \begin{array}{c} \sixvdots \\2\\ \sixvdots  \end{array} \arrow[rd] & &  \\
& & & \begin{array}{c} \vdots \\2\\ \sixvdots  \end{array}\arrow[rd] & &  \begin{array}{c} 1\\1\\1 \end{array} \arrow[rd] & \sixvdots &  \begin{array}{c} 1\\12\\1~ \end{array} \\
&  &  & \sixvdots & \begin{array}{c} 1\\12\\1~ \end{array} \arrow[ru] \arrow[rd] &  & \begin{array}{c} \sixvdots\\ 1\\1 \\ \vdots  \end{array} \arrow[ru] \arrow[rd] & & \\
& & & \begin{array}{c} \sixvdots \\ 1\\1\\ \sixvdots \end{array} \arrow[ru] \arrow[rd] & & \begin{array}{c} 1\\12 \end{array} \arrow[ru] \arrow[rd] & \sixvdots & 1 & \\ 
&&\begin{array}{c} 1\\12 \end{array} \arrow[rd] \arrow[ru] & \sixvdots &   \begin{array}{c} 1 \end{array} \arrow[ru] & &  \begin{array}{c} \sixvdots \\ 1\\2 \\ \vdots \end{array} \arrow[ru] & & & \\
&& &  \begin{array}{c} \sixvdots \\ 1\\2 \\ \sixvdots \end{array} \arrow[ru] & & & & & & & & &
 \end{tikzcd}$.
\end{center} 

%\begin{center}
%$\begin{tikzcd}
%& & 2 \arrow[rd] & & \begin{array}{c} 1\\1\\1\\1 \end{array} \arrow[rd] & & \cdots \\
%\cdots  & \begin{array}{c} 1\\1\\1\\1 \end{array} \arrow[rd] &  & \begin{array}{c} 1\\12\\1~\\1~ \end{array} \arrow[rd]%\arrow[ru] &  &\begin{array}{c} \sixvdots \\ 1\\1\\1 \\ \sixvdots \end{array} & \ldots \\
%& & \begin{array}{c} \vdots \\ 1\\1\\1 \\ \sixvdots \end{array} \arrow[ru] \arrow[rd] & & \begin{array}{c} 1\\12\\1~ \end{array} \arrow[ru] \arrow[rd]&\tenvdots & \cdots\\ 
%\cdots & \begin{array}{c} 1\\1~2\\1~~ \end{array} \arrow[rd]\arrow[ru] & \tenvdots & \begin{array}{c} 1\\1 \end{array} \arrow[rd] \arrow[ru] &  & \begin{array}{c} \sixvdots \\ 1\\1~2 \\ \vdots \end{array} & \cdots \\
%&  &  \begin{array}{c} \sixvdots \\ 1\\1~2 \\ \sixvdots \end{array} \arrow[ru] \arrow[rd] & & 1 \arrow[ru] & & \cdots \\
%\cdots &  1 \arrow[ru] & & \begin{array}{c} 1\\2 \end{array} \arrow[ru] &  & & \\
%\end{tikzcd}$.
%\end{center} 
\noindent
Hence the following facts hold: \\
$(3)$ The indecomposable $\tau$-rigid $B$-modules are the projective
modules $\begin{array}{c} 1\\1~2\\ \vdots ~~ \\ 1~~ \end{array}$ and $2$ and the injective module 
$\begin{array}{c} 1\\ \vdots \\1 \end{array}$ of dimension $n+1$, $1$ and $n$ respectively.\\
$(4)$ The unique non projective $\tau$-tilting $B$-module is the tilting module 
\[
T= \begin{array}{c} 1\\1~2 \\ \vdots~~ \\1~~ \end{array} \oplus \begin{array}{c} 1\\ \vdots \\1  \end{array}\text{
of dimension $2n+1$}. \] 

\noindent
Hence $(i)$ follows from $(1),(2),(3)$ and $(4)$. Finally, the Auslander-Reiten sequence of $A$-modules
$$ 0 \longrightarrow \begin{array}{c} 2\\ \vdots \\2 \end{array} \longrightarrow  
\begin{array}{c} ~~~~~~~2  \\  1 ~~\reflectbox{$\ddots$} ~~~~~ \\ 2~~~~~~~  \end{array} \longrightarrow 1 \longrightarrow 0
$$ 
and the Auslander-Reiten sequence of $B$-modules
$$
0\longrightarrow 2\longrightarrow \begin{array}{c} 1\\1~2\\ \vdots ~~ \\ 1~~ \end{array} 
\longrightarrow \begin{array}{c} 1\\ \vdots \\1 \end{array} \longrightarrow 0
$$
satisfy $(ii)$.

\begin{exam}\label{Ex 7}  
There exist two algebras $A$ and $B\simeq A^{op}$ with the following properties:
\begin{itemize}
\item[(i)] $A$ and $B$ admit indecomposable $\tau$-rigid modules which are 
neither projective, nor injective nor simple.
\item[(ii)] Only one $\tau$-tilting $A$-module is faithful, while any $\tau$-tilting $B$-module 
is faithful.
\item[(iii)] There exists a permutation $\alpha$ (resp. $\beta$) of the indecomposable $\tau$-rigid
$A$-modules (resp. $B$-modules) with the following property:\\

$(+)$ If $X$ and $Y$ are two non isomorphic $\tau$-tilting $A$-modules (resp. $B$-modules), then 
we obtain $Y$ from $X$ by replacing any indecomposable summand $V$ of $X$ which is not a summand of $Y$
by $\alpha(V)$ (resp. $\beta(V)$).    

\item[(iv)] There is a bijection $s$ between the indecomposable $\tau$-rigid $A$-modules and the 
indecomposable $\tau$-rigid $B$-modules such that $\beta = s\alpha s^{-1}$. 
\end{itemize}
\end{exam}
\noindent
{\bf Construction:} Let $A$ be the algebra given by the quiver 
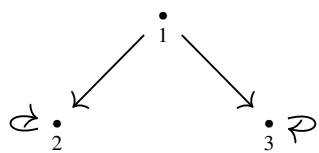
\begin{center}
$\begin{tikzcd}
  &\underset{1}{\bullet} \arrow[dl] \arrow[dr] &  \\
\underset{2}{\bullet}\arrow[loop left] &  & \underset{3}{\bullet}\arrow[loop right]
\end{tikzcd}.$
\end{center} 
with relations $ab=0$ for any paths $a$ and $b$ of length one. Then the Auslander-Reiten quiver of $A$ has the following form:\\

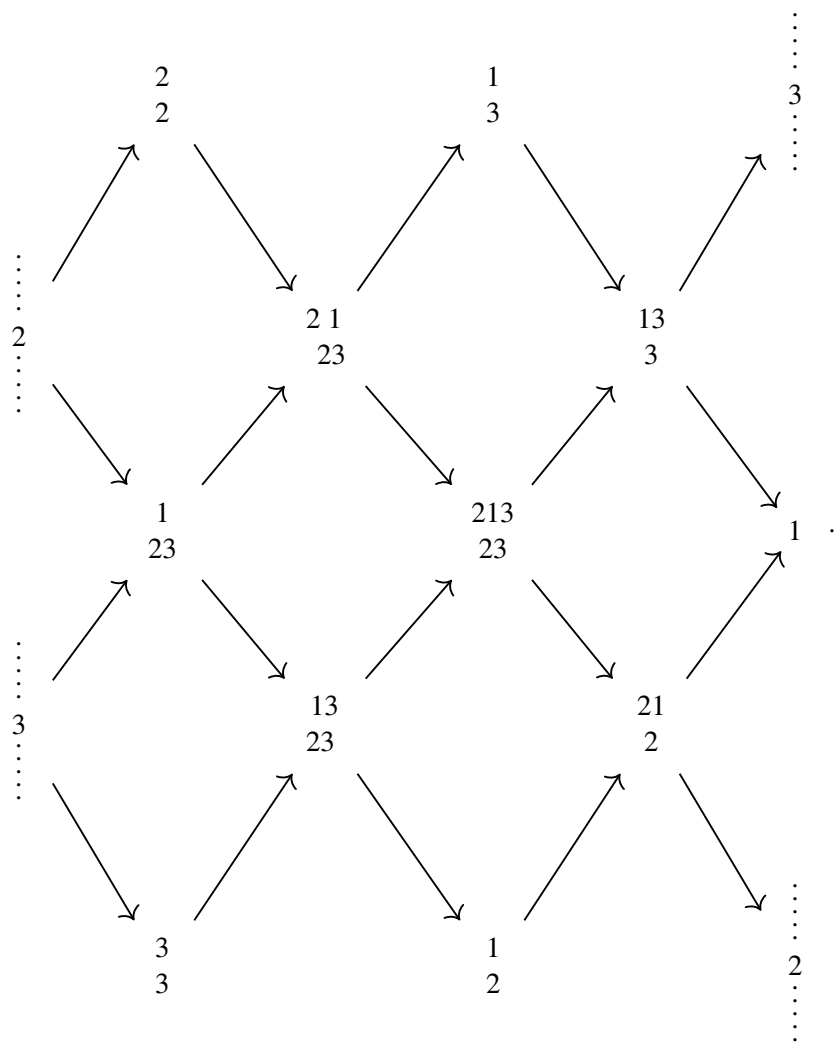
\begin{center}
$\begin{tikzcd}
 & \begin{array}{c} 2\\2 \end{array} \arrow[rd] &  & \begin{array}{c} 1\\3 \end{array} \arrow[rd] &  & \begin{array}{c} \sixvdots \\ 3\\ \sixvdots \end{array}  \\
\begin{array}{c} \sixvdots \\ 2\\ \sixvdots \end{array} \arrow[ru] \arrow[rd]& & \begin{array}{c} ~2~1~\\~~23 \end{array} \arrow[rd] \arrow[ru] &  & \begin{array}{c} 13\\3 \end{array} \arrow[rd] \arrow[ru] & \\
 & \begin{array}{c} 1\\23 \end{array} \arrow[ru] \arrow[rd] & & \begin{array}{c} 213\\~23~ \end{array} \arrow[ru] \arrow[rd] & & 1 \\
\begin{array}{c} \sixvdots \\ 3\\ \sixvdots \end{array} \arrow[ru] \arrow[rd]& & \begin{array}{c} 13\\23~ \end{array}\arrow[ru] \arrow[rd] & & \begin{array}{c} 21\\2 \end{array} \arrow[ru] \arrow[rd] & \\
& \begin{array}{c} 3\\3 \end{array} \arrow[ru] & & \begin{array}{c} 1\\2 \end{array} \arrow[ru] & & \begin{array}{c} \sixvdots \\ 2\\ \sixvdots \end{array}
\end{tikzcd}$.
\end{center} 
Consequently, the following facts hold: \\
\begin{itemize}
\item[(1)] The indecomposable $\tau$-rigid $A$-modules are $\begin{array}{c} 1\\23 \end{array}, \begin{array}{c} 2\\2 \end{array}, \begin{array}{c} 3\\3 \end{array}, \begin{array}{c} 1\\2 \end{array}, \begin{array}{c} 1\\3 \end{array}, 1$.
\item[(2)] The following $A$-modules are not $\tau$-rigid: \\
$$
\begin{array}{c} 1\\23 \end{array} \oplus 1, \begin{array}{c} 2\\2 \end{array} \oplus \begin{array}{c} 1\\3 \end{array},
\begin{array}{c} 2\\2 \end{array} \oplus 1, \begin{array}{c} 3\\3 \end{array} \oplus \begin{array}{c} 1\\2 \end{array},
\begin{array}{c} 3\\3 \end{array} \oplus 1.  
$$ 
\item[(3)] The following $A$-modules are $\tau$-tilting: \\
\[\hspace*{-10mm}
\begin{array}{c} 1\\23 \end{array} \oplus \begin{array}{c} 2\\2 \end{array} \oplus \begin{array}{c} 3\\3 \end{array},
\begin{array}{c} 1\\23 \end{array} \oplus \begin{array}{c} 2\\2 \end{array} \oplus \begin{array}{c} 1\\2 \end{array},
\begin{array}{c} 1\\23 \end{array} \oplus \begin{array}{c} 3\\3 \end{array} \oplus \begin{array}{c} 1\\3 \end{array},
\begin{array}{c} 1\\23 \end{array} \oplus \begin{array}{c} 1\\2 \end{array} \oplus \begin{array}{c} 1\\3 \end{array},
\begin{array}{c} 1\\2 \end{array} \oplus \begin{array}{c} 1\\3 \end{array} \oplus 1.
\]
\end{itemize}
We first deduce from $(1)$ and $(2)$ that the non projective $\tau$-tilting $A$-modules of the form 
$\begin{array}{c} 1\\23 \end{array} \oplus \begin{array}{c} 2\\2 \end{array} \oplus \ldots$ and 
$\begin{array}{c} 1\\23 \end{array} \oplus \begin{array}{c} 3\\3 \end{array} \oplus \ldots$ are the modules 
$\begin{array}{c} 1\\23 \end{array} \oplus \begin{array}{c} 2\\2 \end{array} \oplus \begin{array}{c} 1\\2 \end{array}$
and $\begin{array}{c} 1\\23 \end{array} \oplus \begin{array}{c} 3\\3 \end{array} \oplus \begin{array}{c} 1\\3 \end{array}$,
respectively. 
Next, we observe that the $\tau$-tilting modules of the form $\begin{array}{c} 1\\23 \end{array} \oplus \begin{array}{c} 1\\2 \end{array} \oplus \ldots$(resp. $\begin{array}{c} 1\\23 \end{array} \oplus \begin{array}{c} 1\\3 \end{array} \oplus \ldots$)
are either the module $\begin{array}{c} 1\\23 \end{array} \oplus \begin{array}{c} 1\\2 \end{array} \oplus \begin{array}{c} 2\\2 \end{array}$ (resp. $\begin{array}{c} 1\\23 \end{array} \oplus \begin{array}{c} 1\\3 \end{array} \oplus \begin{array}{c} 3\\3 \end{array}$) or the module $\begin{array}{c} 1\\23 \end{array} \oplus \begin{array}{c} 1\\2 \end{array} \oplus \begin{array}{c} 1\\3 \end{array}$.
Another consequence of $(1)$ and $(2)$ is that the unique $\tau$-tilting of the form $\begin{array}{c} 2\\2 \end{array} \oplus \begin{array}{c} 3\\3 \end{array} \oplus \ldots$ is the projective module $\begin{array}{c} 2\\2 \end{array} \oplus \begin{array}{c} 3\\3 \end{array} \oplus \begin{array}{c} 1\\23 \end{array}$. 
Moreover, the unique $\tau$-tilting modules of the form 
$\begin{array}{c} 1\\2 \end{array} \oplus \begin{array}{c} 1\\3 \end{array} \oplus \ldots$ are $\begin{array}{c} 1\\2 \end{array} \oplus \begin{array}{c} 1\\3 \end{array} \oplus \begin{array}{c} 1\\23 \end{array}$ and $\begin{array}{c} 1\\2 \end{array} \oplus \begin{array}{c} 1\\3 \end{array} \oplus 1$. 
Finally, the unique $\tau$-tilting modules of the form
$\begin{array}{c} x\\x \end{array} \oplus \begin{array}{c} 1\\x \end{array} \oplus \ldots$ and 
$\begin{array}{c} 1\\x \end{array} \oplus 1 \oplus \ldots$ with $x \in \{2,3\}$ are 
$\begin{array}{c} x\\x \end{array} \oplus \begin{array}{c} 1\\x \end{array} \oplus \begin{array}{c} 1\\23 \end{array}$
and $\begin{array}{c} 1\\2 \end{array} \oplus \begin{array}{c} 1\\3 \end{array} \oplus 1$, respectively. 
Putting things together, we conclude that the list in $(3)$ is the list of all $\tau$-tilting $A$-modules.
It follows that \\

\noindent  
$(4)$ $\begin{array}{c} 1\\2 \end{array} $ and $\begin{array}{c} 1\\3 \end{array}$ are uniserial $\tau$-rigid $A$-modules 
which are neither projective, nor injective, nor simple and the projective module 
$\begin{array}{c} 1\\23 \end{array} \oplus \begin{array}{c} 2\\2 \end{array} \oplus \begin{array}{c} 3\\3 \end{array}$
is the unique faithful $\tau$-tilting $A$-module. \\

\noindent
Let $\alpha$ be the following permutation of the indecomposable $\tau$-rigid $A$-modules: 
$$
\begin{array}{c} 1\\23 \end{array} \mapsto 1, 1\mapsto \begin{array}{c} 1\\23 \end{array},
\begin{array}{c} 2\\2 \end{array} \mapsto \begin{array}{c} 1\\3 \end{array}, 
\begin{array}{c} 1\\3 \end{array}\mapsto \begin{array}{c} 2\\2 \end{array},
\begin{array}{c} 3\\3 \end{array} \mapsto \begin{array}{c} 1\\2 \end{array},
\begin{array}{c} 1\\2 \end{array} \mapsto \begin{array}{c} 3\\3 \end{array}.
$$
Since $(3)$ gives the list of all $\tau$-tilting modules, it is to check that \\ 

\noindent
$(5)$ $\alpha$ satisfies property $(+)$ of condition $(iii)$.\\

\noindent
Now, let $B$ be the algebra given by the quiver 
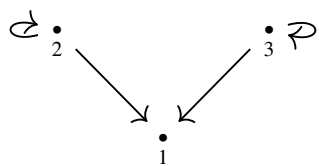
\begin{center}
$\begin{tikzcd}
\underset{2}{\bullet} \arrow[dr] \arrow[loop left]  & & \underset{3}{\bullet} \arrow[dl] \arrow[loop right]  \\
& \underset{1}{\bullet}   & 
\end{tikzcd}.$
\end{center}
with relations $ab=0$ for any paths $a$ and $b$ of length one. 
Then the Auslander-Reiten quiver of $B$ is of the form

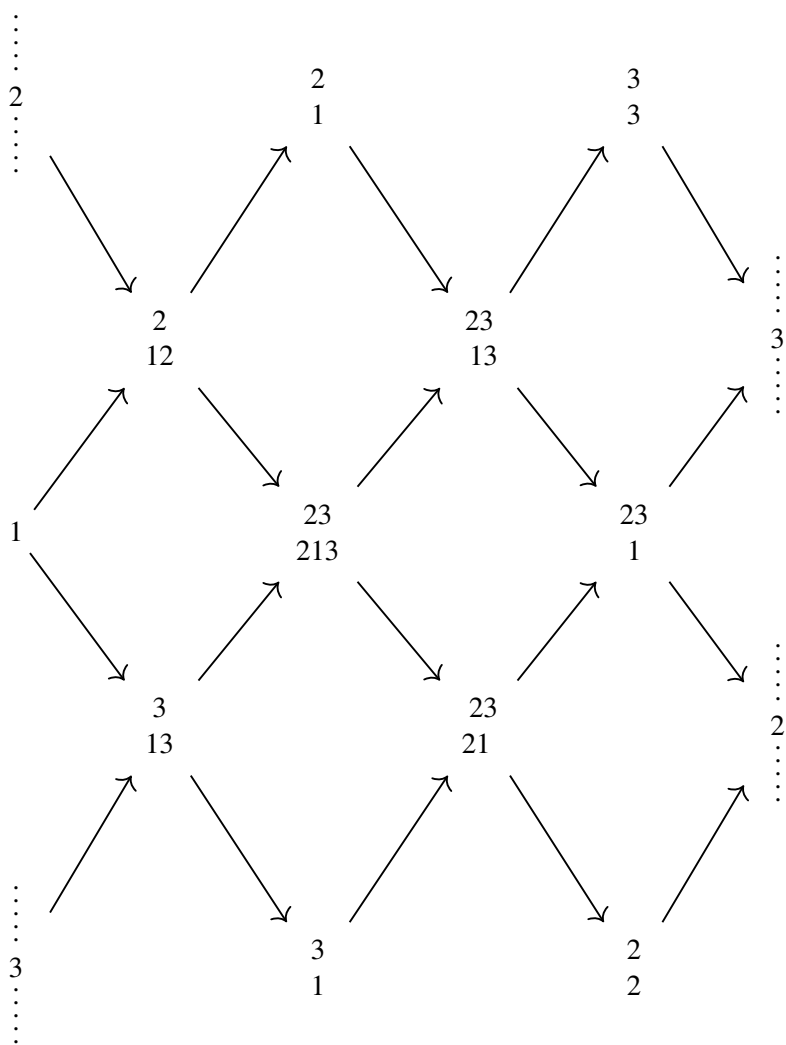
\begin{center}
$\begin{tikzcd}
\begin{array}{c} \sixvdots \\ 2 \\ \sixvdots \end{array} \arrow[rd] &  & \begin{array}{c} 2\\1 \end{array} \arrow[rd] &  & \begin{array}{c} 3\\3 \end{array} \arrow[rd] &    \\
             &  \begin{array}{c} 2\\12 \end{array} \arrow[rd] \arrow[ru] &  & \begin{array}{c} 23\\~13 \end{array} \arrow[rd] \arrow[ru] & & \begin{array}{c} \sixvdots \\ 3 \\ \sixvdots \end{array}\\
1 \arrow[ru] \arrow[rd] & & \begin{array}{c} 23\\213 \end{array} \arrow[ru] \arrow[rd] & & \begin{array}{c} 23\\1 \end{array} \arrow[ru] \arrow[rd] &  \\
 & \begin{array}{c} 3\\13 \end{array}\arrow[rd] \arrow[ru] & & \begin{array}{c} ~23\\21~ \end{array} \arrow[ru] \arrow[rd] & & \begin{array}{c} \sixvdots \\ 2 \\ \sixvdots \end{array} \\
\begin{array}{c} \sixvdots \\ 3 \\ \sixvdots \end{array} \arrow[ru] & & \begin{array}{c} 3\\1 \end{array} \arrow[ru] & & \begin{array}{c} 2\\2 \end{array} \arrow[ru] & 
\end{tikzcd}$.
\end{center}     
Then the following facts hold:\\

\noindent
$(6)$ The indecomposable $\tau$-rigid $B$-modules are: \\
$$
1, \begin{array}{c} 2\\12 \end{array}, \begin{array}{c} 3\\13 \end{array}, \begin{array}{c} 23\\213 \end{array},
\begin{array}{c} 2\\2 \end{array}, \begin{array}{c} 3\\3 \end{array}.
$$ 
$(7)$ The following $B$-modules are not $\tau$-rigid:\\
$$
1\oplus \begin{array}{c} 23\\213 \end{array}, 1\oplus \begin{array}{c} 2\\2 \end{array},
1\oplus \begin{array}{c} 3\\3 \end{array}, \begin{array}{c} 2\\12 \end{array}\oplus \begin{array}{c} 3\\3 \end{array},
\begin{array}{c} 3\\13 \end{array}\oplus \begin{array}{c} 2\\2 \end{array}.  
$$
$(8)$ The following $B$-modules are $\tau$-tilting:\\
\[\scalebox{.75}{$
1\oplus \begin{array}{c} 2\\12 \end{array} \oplus \begin{array}{c} 3\\13 \end{array},
\begin{array}{c} 2\\12 \end{array} \oplus \begin{array}{c} 3\\13 \end{array} \oplus \begin{array}{c} 23\\213 \end{array},
\begin{array}{c} 2\\12 \end{array} \oplus \begin{array}{c} 23\\213 \end{array} \oplus \begin{array}{c} 2\\2 \end{array},
\begin{array}{c} 3\\13 \end{array} \oplus \begin{array}{c} 23\\213 \end{array} \oplus \begin{array}{c} 3\\3 \end{array},
\begin{array}{c} 23\\213 \end{array} \oplus \begin{array}{c} 2\\2 \end{array} \oplus \begin{array}{c} 3\\3 \end{array}.
 $}\]
We first deduce from $(6)$ and $(7)$ that the unique $\tau$-tilting module of the form $1 \oplus \ldots$ is the projective module $1\oplus \begin{array}{c} 2\\12 \end{array} \oplus \begin{array}{c} 3\\13 \end{array}$.
Moreover, the $\tau$-tilting module of the form 
$\begin{array}{c} x\\1x \end{array} \oplus \begin{array}{c} 23\\213 \end{array}\oplus \ldots $ with $x\in \{2, 3\}$
are $\begin{array}{c} x\\1x \end{array} \oplus \begin{array}{c} 23\\213 \end{array}\oplus \begin{array}{c} x\\x \end{array}$
and $\begin{array}{c} 2\\12 \end{array} \oplus \begin{array}{c} 3\\13 \end{array}\oplus \begin{array}{c} 23\\213 \end{array}$.  
We also note that the unique $\tau$-tilting modules of the form 
$\begin{array}{c} 23\\213 \end{array} \oplus \begin{array}{c} x\\x \end{array} \oplus \ldots$ with $x\in \{2, 3\}$
are $\begin{array}{c} 23\\213 \end{array} \oplus \begin{array}{c} x\\x \end{array} \oplus \begin{array}{c} x\\1x \end{array}$
and $\begin{array}{c} 23\\213 \end{array} \oplus \begin{array}{c} 2\\2 \end{array} \oplus \begin{array}{c} 3\\3 \end{array}$.
Next, we deduce from $(6)$ and $(7)$ that the unique non projective $\tau$-tilting module of the form 
$\begin{array}{c} 2\\12 \end{array} \oplus \begin{array}{c} 3\\13 \end{array}\oplus \ldots$ is
$\begin{array}{c} 2\\12 \end{array} \oplus \begin{array}{c} 3\\13 \end{array}\oplus \begin{array}{c} 23\\213 \end{array}$.
On the other hand the unique $\tau$-tilting module of the form 
$\begin{array}{c} 2\\2 \end{array} \oplus \begin{array}{c} 3\\3 \end{array}\oplus \ldots$ is 
$\begin{array}{c} 2\\2 \end{array} \oplus \begin{array}{c} 3\\3 \end{array}\oplus \begin{array}{c} 23\\213 \end{array}$.
We finally note that the unique $\tau$-tilting module of the form 
$\begin{array}{c} x\\1x \end{array} \oplus \begin{array}{c} x\\x \end{array}\oplus \ldots$  with $x \in \{2,3\}$ is 
$\begin{array}{c} x\\1x \end{array} \oplus \begin{array}{c} x\\x \end{array}\oplus \begin{array}{c} 23\\213 \end{array}$.
Putting things together, we conclude that $(8)$ is the list of all  $\tau$-tilting $B$-modules. It follows that \\

\noindent
$(9)$ $\begin{array}{c} 23\\213 \end{array}$ is the unique indecomposable $\tau$-rigid $B$-module which is neither projective, nor injective, nor simple and any $\tau$-tilting $B$-module is faithful. \\
\\
Now, let $\beta$ denote the following permutation of the indecomposable 
$\tau$-rigid $B$-modules: \\
$$
1\mapsto \begin{array}{c} 23\\213 \end{array}, \begin{array}{c} 23\\213 \end{array} \mapsto 1,
\begin{array}{c} 2\\21 \end{array} \mapsto \begin{array}{c} 3\\3 \end{array},
\begin{array}{c} 3\\3 \end{array} \mapsto \begin{array}{c} 2\\21 \end{array},
\begin{array}{c} 3\\31 \end{array} \mapsto \begin{array}{c} 2\\2 \end{array},
\begin{array}{c} 2\\2 \end{array} \mapsto \begin{array}{c} 3\\31 \end{array}.
$$ 
Then it is easy to check that \\
\\
$(10)$ $\beta$ satisfies property $(+)$ of condition $(iii)$.\\
\\
Hence $(i)$ and $(ii)$ follow from $(4)$ and $(9)$, while 
$(iii)$ follows from $(5)$ and $(10)$. 
Finally, let $s$ denote the following bijection between the indecomposable 
$\tau$-rigid $A$-modules and the indecomposable $\tau$-rigid $B$-modules: \\
$$
\begin{array}{c} 1\\23 \end{array} \mapsto \begin{array}{c} 23\\213 \end{array},
\begin{array}{c} 2\\2 \end{array} \mapsto \begin{array}{c} 2\\2 \end{array},
\begin{array}{c} 3\\3 \end{array} \mapsto \begin{array}{c} 3\\3 \end{array},
\begin{array}{c} 1\\2 \end{array} \mapsto \begin{array}{c} 2\\12 \end{array},
\begin{array}{c} 1\\3 \end{array} \mapsto \begin{array}{c} 3\\13 \end{array},
1 \mapsto 1.
$$
Then we have $\beta = s\alpha s^{-1}$, and so $(iv)$ holds. 

\begin{rem}
\normalfont{
Looking at the two Auslander-Reiten quivers in Example \ref{Ex 7}, we see that the indecomposable $\tau$-rigid
$A$-modules are $P_1$, $P_2$, $P_3$, $I_1 =S_1 \simeq P_1/ (S_2\oplus S_3)$, $\begin{array}{c} 1\\2 \end{array}\simeq P_1/S_3$ and $\begin{array}{c} 1\\3 \end{array} \simeq P_1/S_2$.\\
\\
On the other hand, the indecomposable $\tau$-rigid $B$-modules are $P_1=S_1$, $P_2$, $P_3$, 
$\begin{array}{c} 23\\213 \end{array} \simeq (P_2\oplus P_3)/S_1$, 
$I_2=\begin{array}{c} 2\\2 \end{array}\simeq P_2/S_1$ and $I_3 = \begin{array}{c} 3\\3 \end{array} \simeq P_3/S_1$.
Hence, we may roughly speaking say that we obtain the indecomposable non-projective $\tau$-rigid $B$-modules from the indecomposable non-projective $\tau$-rigid $A$-modules by interchanging "numerators" with "denominators" and projective modules with simple modules.}  
\end{rem}

\begin{rem}
\normalfont{
The bijection $s$ satisfying the last condition of Example \ref{Ex 7} induces a bijection 
between $\tau$-tilting $A$-modules and $\tau$-tilting $B$-modules.
If $V$ is an indecomposable $\tau$-rigid $A$-module, then the following 
facts hold:\\

$\bullet$ $s(V)$ is projective $\Longleftrightarrow$ $V$ is not projective. \\

$\bullet$ $s(V)$ is simple $\Longleftrightarrow$ $V$ is simple.\\

$\bullet$ $s(V)$ is sincere $\Longleftrightarrow$ $V$ is sincere.\\

$\bullet$ $s(V)$ is uniserial $\Longrightarrow$ $V$ is uniserial. \\

$\bullet$ $s(V)$ is faithful $\nRightarrow$ $V$ is faithful. \\

$\bullet$ $dim_{K}(V) \leq dim_{K}s(V) \leq dim_{K}(V)+2$.}
\end{rem}

%\begin{exam}\label{Ex 11}
%There is an algebra $A$ admitting a uniserial $\tau$-rigid module $U$ with the following properties:
%\begin{itemize}
%\item[(i)] The $\tau$-orbit of $U$ contains of $\tau$-rigid modules.
%\item[(ii)] $U$ has infinite projective and injective dimension.
%\item[(iii)] $Ext^{n}_{A}(U,U)\neq 0$ for any $n\geq 2$.
%\item[(iv)] The Auslander-Reiten quiver $A$ contains two connected components $\mathcal{C}$ and $\mathcal{D}$ 
%with exactly two indecomposable $\tau$-rigid modules. Moreover, there exist exactly three $\tau$-tilting modules
%$T_1, T_2, T_3$ with indecomposable summands in $\mathcal{C} \cup \mathcal{D}$.
%\item[(v)] 
%\end{itemize}
%\end{exam} 

A remark of \cite[Section 4.1]{IR} points out that a selfinjective algebra, hence an algebra with exactly one basic tilting module, "usually has a lot of $\tau$-tilting modules". The next result shows that the dimension of a $\tau$-tilting module may be rather small with respect to that of the unique basic tilting module.

\begin{exam}\label{Ex 10} 
There exist a selfinjective $K$-algebra $A$ with six simple modules and $A$-modules $T$ and $T'$ with the following properties:
\item[(i)] Every indecomposable projective $A$-module has dimension 7.
\item[(ii)] $T$ and $T'$ are basic $\tau$-tilting modules of infinite projective dimension 
of the form $1\oplus3\oplus5\oplus...$ and $\begin{array}{c} 1\\2 \end{array} \oplus \begin{array}{c} 4\\5 \end{array}\oplus ...$ with exactly one indecomposable projective summand,  and the dimension of the indecomposable non projective summands of $T$ and $T'$ are $1,1,1,3,5$ and $1,1,2,2,5$ respectively.
\item[(iii)] $T$ and $T'$ have a maximal common summand $U$ of dimension 13.
\item[(iv)] If $d=1$ (resp. $d=2$) and $V$ is a basic $\tau$-rigid module which is the direct sum of indecomposable modules of dimension $d$, then we have $dim_{K}(V)\leq 3$ (resp. $dim_{K}(V)\leq 4$).
\item[(v)] If $d=3,4,5$, then the direct sum of two non isomorphic indecomposable modules of dimension $d$ is not $\tau$-rigid.
\item[(vi)] $T$ and $T'$ generate the same number of indecomposable modules. 
\item[(vii)] $T$(resp. $T'$) is generated by three (resp. two) indecomposable modules.
\item[(viii)] There is a bijection $s$ between the indecomposable summands $T_1, \ldots, T_6$ of $T$ and the indecomposable summands of $T'$ such that either $s(T_{i})\simeq T_{i}$ or \\$Ext_{A}^1(T_{i}, s(T_{i}))\oplus Ext_{A}^{1}(s(T_{i}), T_{i}) \neq 0$ for any $i$.
\item[(ix)] If $X$ is a $\tau$-tilting module of the form $1\oplus 3\oplus 5\oplus ...$ or $\begin{array}{c} 1\\2 \end{array} \oplus \begin{array}{c} 4\\5 \end{array}\oplus ...$, then we have $dim_{K}(X)\geq 18=dim_{K}(T)=dim_{K}(T')$.
\end{exam}

\noindent
{\bf Construction:} Let $A$ be the algebra given by the quiver

$$
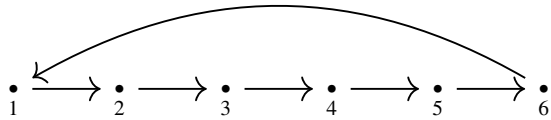
\begin{tikzcd}
\underset{1}{\bullet} \arrow[r] & \underset{2}{\bullet} \arrow[r] & \underset{3}{\bullet} \arrow[r] &
\underset{4}{\bullet} \arrow[r] & \underset{5}{\bullet} \arrow[r] & \underset{6}{\bullet} \arrow[lllll, bend right=30]
\end{tikzcd}$$
with the property that the composition of seven arrows is always zero. Then for any vertex $x$, the indecomposable projective module $P_{x}$ is the uniserial module of the form $\begin{array}{c} x\\\vdots\\x  \end{array}$ of the dimension $7$.
Hence, $(i)$ holds. Moreover, the shape of the Auslander-Reiten quiver of $A$ shows that \\
\\
$(1)$ The indecomposable $\tau$-rigid modules are the indecomposable modules of dimension different from six.\\

\noindent
Let $T$ and $T'$ denote the following modules: 
\begin{center}
$T= P_1 \, \oplus \, 1 \oplus \begin{array}{c} 1\\2\\3 \end{array} \oplus \, 3 \, \oplus \, \begin{array}{c} 1\\2\\3\\4\\5 \end{array} \oplus \, 5$,  
$T'= P_1 \oplus \, 1 \oplus \begin{array}{c} 1\\2 \end{array} \oplus \begin{array}{c} 1\\2\\3\\4\\5 \end{array} \oplus \, 4 \, \oplus \begin{array}{c} 4\\5 \end{array}$.
\end{center}

Then $T$ and $T'$ are $\tau$-tilting modules. Since the simple module $1$ has infinite projective dimension, it follows that 
$(ii)$ holds.
We next observe that the module $U=P_1\oplus 1\oplus \begin{array}{c} 1\\ \vdots \\5 \end{array}$ satisfies $(iii)$. Moreover, $1\oplus 3 \oplus 5$ and $2 \oplus 4\oplus 6$ are the basic semisimple $\tau$-rigid modules with the largest possible dimension.
On the other hand, $\begin{array}{c} 1\\2 \end{array} \oplus \begin{array}{c} 4\\5 \end{array}$, $\begin{array}{c} 2\\3 \end{array} \oplus \begin{array}{c} 5\\6 \end{array}$ and $\begin{array}{c} 3\\4 \end{array} \oplus \begin{array}{c} 6\\1 \end{array}$  are the basic $\tau$-rigid modules which are direct sums of indecomposable modules of dimension two and have the largest possible dimension. Hence, $(iv)$ holds and a direct calculation also shows that also $(v)$ holds. Finally, the indecomposable quotients of $P_1$ and the modules $3$ and $5$ (resp.$\begin{array}{c} 4\\5 \end{array}$ and $4$) are the indecomposable modules generated by $T$(resp. $T'$). Hence $(vi)$ holds.
Since the modules $P_1\oplus 3\oplus 5$ and $P_1\oplus \begin{array}{c} 4\\5 \end{array}$ generate $T$ and $T'$ respectively, it follows that $(vii)$ holds. Let $s$ be the bijective map such that $P_1\mapsto P_1$, $1\mapsto 1$, $\begin{array}{c} 1\\ \vdots \\5 \end{array} \mapsto \begin{array}{c} 1\\ \vdots \\5 \end{array}$, $3\mapsto \begin{array}{c} 1\\2 \end{array}$, $5\mapsto 4$ and $\begin{array}{c} 1\\2\\3 \end{array} \mapsto \begin{array}{c} 4\\5 \end{array}$. Then $s$ obviously satisfies $(viii)$. Let $X= 1\oplus 3\oplus 5\oplus ...$ be a $\tau$-tilting module. Next, let $W$ be an indecomposable module of dimension $d\in \{2,4\}$. Then it is easy to see that $1\oplus 3\oplus 5\oplus W$ is not $\tau$-rigid. This observation and $(v)$ imply that \\
\\
$(2)$ $X= 1\oplus 3\oplus 5\oplus ...$ has at most one indecomposable summand of dimension $3$ or $5$, but $X$ doesn't have indecomposable summands of dimension $2$ or $4$. \\
\\
Putting $(1)$ and $(2)$ together, we deduce from $(iv)$ that\\
\\
$(3)$ $dim_{K}(X)=18 \geq dim_{K}(T)=dim_{K}(T')$. \\
\\
Finally, let $X=\begin{array}{c} 1\\2 \end{array}\oplus \begin{array}{c} 4\\5 \end{array}\oplus ...$ be a basic $\tau$-tilting module. Then it is easy to see that $X$ does not have indecomposable summands of dimension $3$ and $4$. Moreover, we already know that $X$ does not contain a semisimple summand of the form $a\oplus b\oplus c$.
Consequently, $dim_{K}(X)$ is as small as possible if $X=\begin{array}{c} 1\\2 \end{array}\oplus \begin{array}{c} 4\\5 \end{array}\oplus W$, where the indecomposable summands of $W$ have dimension $1,1,5,7$. Hence, we have \\
\\
$(4)$ $dim_{K}(X)\geq 18=dim_{K}(T)=dim_{K}(T')$.\\
\\
Hence, $(ix)$ follows from $(3)$ and $(4)$.

 %\begin{tikzpicture}

  %  \underset{1}{\node} (1) at (0,0){1};

   % \node (2) at (1,1){2};

    %\node (3) at (1,0){3};
     
    %\node (4) at (1,-1){4};
    
    %\draw[->] (1) -- (2);
    %\draw[->] (1) -- (3);
     %\draw[->] (1) -- (4);

 %\end{tikzpicture}
\begin{rem}
\normalfont{
The $\tau$-tilting modules $T$ and $T'$ of Example \ref{Ex 10} have many indecomposable direct summands of small dimension. 
Indeed, the direct sum of the indecomposable direct summands of $T$ (resp. $T'$) of dimension $\leq 2$ has dimension $3$ (resp.  $6$). It is easy to find $\tau$-tilting modules $W$ and $W'$ of dimension $19$ and $20$ respectively with exactly two simple
summands and one indecomposable summand of dimension two. Indeed, the modules  
$$W = 1 \oplus 4 \oplus \begin{array}{c} 1\\2 \end{array} \oplus \begin{array}{c} 6\\1\\2 \end{array} \oplus \begin{array}{c} 6\\1\\2\\3\\4 \end{array} \oplus P_6 \text{ and } W' = 1 \oplus 4 \oplus \begin{array}{c} 1\\2 \end{array} \oplus \begin{array}{c} 1\\2\\3\\4 \end{array} \oplus \begin{array}{c} 6\\1\\2\\3\\4 \end{array} \oplus P_6$$  have all the desired properties.}
\end{rem}

\textbf{Acknowledgements:}

We wish to thank the organizers of the following conferences on Representation Theory:
\begin{itemize}
\item VirtARTA 2021, where the second author gave an online talk on \cite{GM};
\item ICRA 2022, where the first author gave an online talk on this paper.
\end{itemize}
We wish to thank also
\begin{itemize}
\item Bernhard Keller, who asked "Can you replace tilting modules by $\tau$-tilting modules?" in the pages of remarks concerning
the presentation at VirtARTA 2021.
\item Ramin Ebrahimi, who asked "Can you replace Ext-rigidity by $\tau$-rigidity?" after the presentation at ICRA 2022 and later
suggested us to look at \cite[Theorem 2.12]{AIR}.
\end{itemize}
Keller and Ebrahimi's questions have a positive answer. Indeed, we cannot simply replace "tilting modules" by "$\tau$-tilting modules" in \cite[Theorem 1]{GM}. However, if we make the following substitutions
\begin{itemize}
\item tilting module $\mapsto$ $\tau$-tilting module
\item partial tilting module $\mapsto$ $\tau$-rigid module
\item $Ext^1_{A}(X_{i},Y_{s(i)})\oplus  Ext^1_{A}(Y_{s(i)},X_{i})\neq 0 \mapsto \\ \phantom{.}\hfill Hom_{A}(X_{i},\tau Y_{s(i)})\oplus 
Hom_{A}(Y_{s(i)},\tau X_{i})\neq 0$
\end{itemize}
then \cite[Theorem 2.12]{AIR} is the $\tau$-tilting version of Ringel's result \cite[page 166]{R}, which says the following: If $T$ is a tilting-module and $T\oplus S$ is a partial tilting module, then $S \in addT$.
Consequently, after a word by word repetition of the proof of \cite[Theorem 2.5]{GM}, we obtain the proof of the following result: \\
Given two $\tau$-tilting $A$-modules $X$ and $Y$ of the form $X= X_1\oplus \ldots \oplus X_{n}$ and  $Y = Y_1\oplus \ldots \oplus Y_{n}$, with $X_{i}$ and $Y_{i}$ indecomposable for any $i$, there exists a permutation $s\in S_{n}$ such that either
$X_{i}\simeq Y_{s(i)}$ or $Hom_{A}(X_{i},\tau Y_{s(i)})\oplus Hom_{A}(Y_{s(i)},\tau X_{i})\neq 0$.

\end{document}